\documentclass[12pt]{amsart}

\usepackage{amssymb}
\usepackage{latexsym}
\usepackage{amsmath}
\usepackage{paralist}
\usepackage[all]{xy}

\newtheorem{theorem}{Theorem}[section]

\newtheorem{corollary}[theorem]{Corollary}

\newtheorem{definition}[theorem]{Definition}
\newtheorem{example}[theorem]{Example}

\newtheorem{applications}[theorem]{Applications}

\newtheorem{lemma}[theorem]{Lemma}

\newtheorem{proposition}[theorem]{Proposition}
\newtheorem{remark}[theorem]{Remark}

\newcommand{\bbC}{{\mathbb{C}}}

\newcommand{\bbR}{{\mathbb{R}}}
\newcommand{\bbT}{{\mathbb{T}}}
\newcommand{\bbZ}{{\mathbb{Z}}}

  \newcommand{\A}{{\mathcal{A}}}
  \newcommand{\B}{{\mathcal{B}}}
  \newcommand{\C}{{\mathcal{C}}}
  \newcommand{\D}{{\mathcal{D}}}
  
  \newcommand{\F}{{\mathcal{F}}}
  
\renewcommand{\H}{{\mathcal{H}}}

  \newcommand{\M}{{\mathcal{M}}}

\renewcommand{\S}{{\mathcal{S}}}
  \newcommand{\T}{{\mathcal{T}}}

\newcommand{\fA}{{\mathfrak{A}}}

\newcommand{\fC}{{\mathfrak{C}}}

\newcommand{\fT}{{\mathfrak{T}}}

\def\gs{\sigma}
\def\ga{\alpha}
\def\gl{\lambda}
\def\om{\omega}
\def\gd{\delta}
\def\gi{\iota}

\newcommand{\sca}[1]{\left\langle#1\right\rangle}

\newcommand\Ref{\mathop{\rm Ref}}

\begin{document}

\title[Semicrossed Products and Reflexivity]{Semicrossed Products and Reflexivity}

\author{Evgenios T.A. Kakariadis}

\address{Evgenios T.A. Kakariadis\\ Department of Mathematics,\\
University of Athens, \\
Panepistimioupolis, \\GR-157 84, Athens, GREECE}
   \email{mavro@math.uoa.gr}

\thanks{\ The author was supported by an SSF scholarship.}


\date{}

\maketitle

\begin{abstract}
Given a w*-closed unital algebra $\A$ acting on $H_0$ and a
contractive w*-continuous endomorphism $\beta$ of $\A$, there is a
w*-closed (non-selfadjoint) unital algebra
$\mathbb{Z}_+\overline{\times}_\beta \A$ acting on
$H_0\otimes\ell^2({\mathbb{Z}_+})$, called the w*-semicrossed
product of $\A$ with $\beta$. We prove that
$\mathbb{Z}_+\overline{\times}_\beta \A$ is a reflexive operator
algebra provided $\A$ is reflexive and $\beta$ is unitarily
implemented, and that $\mathbb{Z}_+\overline{\times}_\beta \A$ has
the bicommutant property if and only if so does $\A$. Also, we show
that the w*-semicrossed product generated by a commutative
C*-algebra and a *-endomorphism is reflexive.
\end{abstract}

\section*{Introduction}

\noindent

As is well known, to construct the
C*-crossed product of a unital C*-algebra $\C$ by a *-isomorphism
$\ga: \C\rightarrow \C$, we begin with the Banach space
$\ell^1(\bbZ,\C,\ga)$ which is the closed linear span of the
monomials $\gd_n\otimes x$, $n\in\bbZ$, $x\in\C$, under the norm
$|\sum_{n=-k}^{k}\gd_n\otimes x_n|_1= \sum_{n=-k}^{k} \|x_n\|_\C$,
equipped with the (isometric) involution $(\gd_n\otimes
x)^*=\gd_{-n}\otimes \ga^{-n}(x^*)$. Now, there are two ``natural''
ways  to define multiplication in $\ell^1(\bbZ,\C,\ga)$; either the
left multiplication $
 (\gd_n\otimes x) \ast_l (\gd_m\otimes y) = \gd_{n+m}\otimes a^m(x)y
$, or the right one $
 (\gd_n\otimes x) \ast_r  (\gd_m\otimes y) =\gd_{n+m}\otimes xa^n(y)
$. Then the corresponding algebras are isometrically *-isomorphic
via the map $\Psi(\gd_n\otimes x)=\gd_{-n}\otimes a^{-n}(x)$. We can
see that $(\ell^1(\bbZ,\C,\ga)_l)^{opp}
=\ell^1(\bbZ,\C^{opp},\ga)_r$, where for an algebra $\B$, $\B^{opp}$
is the space $\B$ along with the multiplication $x\odot y:=yx$;
hence, in case $\C$ is commutative, each algebra is the opposite of
the other. The left and right crossed product are the completion of
the corresponding involutive Banach algebras under a universal norm
induced by the $|\cdot|_1$-contractive *-representations (hence,
they are C*-algebras characterized by a universal property) and the
map $\Psi$ extends to a C*-isomorphism. Moreover, it can be proved
that the crossed product is *-isomorphic to the reduced crossed
product $C^*_l(\C)$, i.e. the norm closure of the range of the left
regular representation, and thus we end up with just one object to
which we refer as \emph{the crossed product of the dynamical system}
$(\C,\ga)$. The key fact is that there is a bijection between the
$|\cdot|_1$-contractive *-representations of each of these
$\ell^1$-algebras and the (left or right) covariant unitary pairs
(see section 1).

If we wish to construct a non-selfadjoint analogue, we can see that
there are more possibilities. For example, Peters defined the
semicrossed product as the completion of the Banach algebra
$\ell^1(\bbZ_+,\C,\ga)_l$ under the universal norm that arises from
the left covariant isometric pairs and examined the case when $\ga$
is an injective *-endomorphism of $\C$. He proved that this
semicrossed product embeds isometrically in a crossed product (see
\cite{Pet84}) and, for the commutative case, that this crossed
product is the C*-envelope of the semicrossed product (see
\cite{Petarx}).

In section 1 we use an alternative definition using ``sufficiently
many'' homomorphisms of the Banach algebra $\ell^1(\bbZ_+,\C,\ga)_l$
(see also \cite{Dav07}). The advantage is that there is a bijection
between the left covariant contractive pairs and the homomorphisms
of the Banach algebra $\ell^1(\bbZ_+,\C,\ga)_l$. Moreover, there is
a duality between the left covariant contractive pairs and the right
covariant contractive pairs, which induce the homomorphisms of the
Banach algebra $\ell^1(\bbZ_+,\C,\ga)_r$; hence, we get similar
results for the right version. Also, using a dilation theorem of
\cite{MuhSol06}, we can see that this definition gives the one in
\cite{Pet84}. If we consider the maximal operator space structure,
then the semicrossed products are operator algebras with a universal
property that characterizes them up to completely isometric
isomorphism. In theorem \ref{onlyone} we prove that the semicrossed
product is independent of the way $\C$ is (faithfully) represented
and in theorem 1.5 we prove that in case $\ga$ is a *-isomorphism,
its C*-envelope is exactly the crossed product. So, in order to
define a w*-analogue of the semicrossed product that arises by a
w*-continuous contractive endomorphism $\beta$ of a w*-closed
subalgebra $\A$ of some $\B(H_0)$ (for example, a von Neumann
algebra), either we take the w*-closed linear span of a
non-selfadjoint left regular representation or the w*-closed linear
span of the analytic polynomials of the von Neumann crossed product,
depending on the properties of $\beta$.

In section 2 we analyze the properties of the w*-semicrossed
product, in case $\beta$ is unitarily implemented. First of all, we
study the connection between the semicrossed product and the
w*-tensor product $\A\overline{\otimes}\T$, where $\T$ is the
algebra of the analytic Toeplitz operators, and give an example when
these two algebras are incomparable. A main result of this section
is the reflexivity of the w*-semicrossed product, when $\A$ is
reflexive. Recall that a subspace $\S\subseteq B(H)$ is
\emph{reflexive} if it coincides with its \emph{reflexive cover},
namely $\Ref(\S)=\{T\in\B(H): T\xi\in \overline{\S\xi}, \text{for all}\
\xi\in\H\}$ (see \cite{LogSu75}); unlike \cite{LogSu75},
we will call $\S$ \emph{hereditarily reflexive}
if every w*-closed subspace of $\S$ is reflexive.  As a consequence
we have that, when a unitary implementation condition holds, the
w*-closed image of $lt_\pi$ induced by a representation $(H_0,\pi)$
of $\C$ is reflexive. Also, we get several known results as
applications. As another main result, we prove that the
w*-semicrossed product is the commutant of a w*-semicrossed product
and is its own bicommutant if and only if the same holds for $\A$.

In the last section we consider the semicrossed product of a
commutative C*-algebra $C(K)$ with a continuous map
$\phi:K\rightarrow K$. As observed in theorem \ref{onlyone}, the
representations induced by a character of $C(K)$, say $ev_t$, $t\in
K$, suffice to obtain the norm of the semicrossed product and play a
significant role for its study. First, we show that the w*-closure
of such representations is always reflexive; in fact, it has the
form $(\fT P_{n_0})\oplus (\T P_0) \oplus \cdots \oplus (\T
P_{p-1})$, where $\fT$ is the algebra of lower triangular operators
in $\B(\ell^2(\bbZ_+))$, $\T$ is the algebra of analytic Toeplitz
operators and $P_{n_0}, P_0, \dots P_{p-1}$ some projections
determined by the orbit of the point $t\in K$.\\

In what follows we use standard notation, as in \cite{KadRin86}
for example. $\bbZ_+=\{0,1,2,\dots\}$ and all infinite sums are
considered in the strong-convergent sense. Throughout, we use the
symbol $v$ for the unilateral shift on $\B(\ell^2(\bbZ_+))$, given
by $v(e_n)=e_{n+1}$. A useful tool for the proofs in sections 2 and
3 is a \emph{F\'{e}jer-type Lemma}; consider the unitary action of
$\bbT$ on $H=H_0\overline{\otimes}\ell^2(\bbZ_+)$ induced by the
operators $U_s,\ s\in \bbR$, given by $U_s(\xi\otimes
e_n)=e^{ins}\xi\otimes e_n$. For every $T\in \B(H)$ and every $m\in
\bbZ$ we define the \emph{``m-Fourier coefficient''}
 \[
 G_m(T)=\int_0^{2\pi} U_sTU_s^*e^{-ims}\frac{ds}{2\pi},
 \]
the integral taken as the w*-limit of Riemann sums. If we set
$\sigma_l(T)(t) = \frac{1}{l+1} \sum_{n=0}^l \sum_{m=-n}^{n} G_m(T)
\exp(imt)$, then $\gs_l(T)(0)\stackrel{w*}{\to} T$. Note that
$G_m(\cdot)$ is w*-continuous for every $m\in\bbZ$.

Now, for every $\kappa, \gl \in \bbZ_+$, and $T\in\B(H)$ let the
\emph{``matrix elements''} $T_{\kappa,\gl}\in \B(H_0)$ be defined by
 $\sca{T_{\kappa,\gl}\xi,\eta}=\sca{T(\xi\otimes e_\gl),\eta \otimes
 e_\kappa}$, $\xi,\eta\in H_0$; then we can write the Fourier
coefficients explicitly by the formula
\begin{align*}
G_m(T)=
\begin{cases}
   V^m(\sum_{n\geq 0} T_{m+n,n}\otimes p_n) &,\ \text{when}\ m\geq 0,\\
   (\sum_{n\geq 0} T_{n,-m+n}\otimes p_n) (V^*)^ {-m} &,\ \text{when}\
   m<0,
\end{cases}
\end{align*}
where $V=1_{H_0}\otimes v$. For simplicity, we define the diagonal
matrices
\begin{align*}
T_{(m)}=
\begin{cases}
   \sum_{n\geq 0} T_{m+n,n}\otimes p_n &,\ \text{when}\ m\geq 0,\\
   \sum_{n\geq 0} T_{n,-m+n}\otimes p_n &,\ \text{when}\ m<0,
\end{cases}
\end{align*}
Note that the sums converge in the w*-topology as well, since the
partial sums are uniformly bounded by $\|T\|$. Hence, $G_m(T)$ is
the \emph{m-diagonal} of $T$, when we view $H$ as the $\ell^2$-sum
of copies of $H_0$.

\section{Semicrossed products of C*-algebras}

\noindent

Let $\C$ be a unital C*-algebra and $\ga: \C\rightarrow \C$ a
*-morphism; define $\ell^1(\bbZ_+,\C,\ga)$ to be the closed linear
span of the monomials $\gd_n\otimes x$, $n\in\bbZ_+$, $x\in\C$,
under the norm
\[
 \left|\sum_{n=0}^{k}\gd_n\otimes x_n\right|_1= \sum_{n=0}^{k}
 \|x_n\|_\C.
\]
We endow $\ell^1(\bbZ_+,\C,\ga)$ either with the left multiplication
$
 (\gd_n\otimes x) \ast_l (\gd_m\otimes y) = \gd_{n+m}\otimes a^m(x)y,
$ or with the right one $
 (\gd_n\otimes x) \ast_r  (\gd_m\otimes y) =\gd_{n+m}\otimes
 xa^n(y),$ and denote the corresponding Banach algebras by
$\ell^1(\bbZ_+,\C,\ga)_l$ and $\ell^1(\bbZ_+,\C,\ga)_r$,
respectively. One can see that $(\ell^1(\bbZ_+,\C,\ga)_l)^{opp}$ is
exactly $\ell^1(\bbZ_+,\C^{opp},\ga)_r$, where, if $\B$ is an
algebra, $\B^{opp}$ is the space $\B$ with the multiplication
$x\odot y:=yx$. Thus, in case $\C$ is commutative,
each algebra is the opposite of the other.\\

Let $(H,\pi)$ be a *-representation of $\C$ and $T$ a contraction in
$\B(H)$. The pair $(\pi,T)$ is called a \emph{left covariant
contractive (l-cov.con.)} pair, if the left covariance relation is
satisfied, i.e. $\pi(x)T=T\pi(\ga(x)),\ x\in\C$. If, in particular,
$T$ is an isometry, pure isometry, co-isometry or unitary, then we
will call such a pair \emph{a left covariant isometric, purely
isometric, co-isometric or unitary pair}. We can see that every
l-cov.con. pair induces a contractive representation
$(H,T\times\pi)$ of $\ell^1(\bbZ_+,\C,\ga)_l$, given by
 \[
 (T\times\pi)\left(\sum_{n=0}^k \gd_n\otimes x_n\right)=
 \sum_{n=0}^k T^n\pi(x_n).
 \]
Conversely, if $\rho:\ell^1(\bbZ_+,\C,\ga)_l\rightarrow \B(H)$ is a
contractive representation, then $(H,\rho)$ restricts to a
contractive representation $(H,\pi)$ of the C*-algebra $\C$, thus a
*-representation. If we set $\rho(\gd_1\otimes e)=T$, then
$\|T^n\|\leq 1$, for every $n\in\bbZ_+$. It is easy to check that
the pair $(\pi,T)$ satisfies the left covariance relation.

Analogously, there is a bijection between \emph{right covariant contractive (r-cov.con.)}
pairs $(\pi, T)$ 
(i.e. satisfying the right covariance condition
$T\pi(x)=\pi(\ga(x))T,\ x\in\C$) and  contractive representations
$\pi\times T$ of the algebra $\ell^1(\bbZ_+,\C,\ga)_r$. Note that if
$(\pi,T)$ is a l-cov.con. pair then $(\pi,T^*)$ is a r-cov.con.
pair. Thus $TT^*$ commutes with $\pi(\C)$.

\begin{example}
\textup{Let $(H_0,\pi)$ be a faithful *-representation of $\C$ and
define on $H_0\otimes\ell^2(\bbZ_+)$ the representation
$\widetilde{\pi}(x)=\text{diag}\{\pi(\ga^n(x)):n\in\bbZ_+\}$ and
$V=1_{H_0}\otimes v$, where $v$ is the unilateral shift. Then
$(\widetilde{\pi},V)$ is a l-cov.is. pair. For simplicity we will
denote the corresponding representation $V\times\widetilde{\pi}$, by
$lt_\pi$. As mentioned before, the pair $(\widetilde{\pi}, V^*)$ is
a r-cov.con. pair which induces the representation
$rt_\pi:=\widetilde{\pi}\times V^*$. One can check that $lt_\pi$ and
$rt_\pi$ are faithful.}
\end{example}

\begin{definition}
The (left) semicrossed product $\bbZ_+\times_\ga\C$ is the
completion of $\ell^1(\bbZ_+,\C,\ga)_l$ under the norm
\[
 \|F\|_l=\sup\{\|(T\times\pi)(F)\|: (\pi,T)\ \textup{is a l-cov.con.
 pair}\}.
\]
Analogously, the (right) semicrossed product $\C\times_\ga\bbZ_+$ is
the completion of $\ell^1(\bbZ_+,\C,\ga)_r$ under the norm
\[
 \|F\|_r=\sup\{\|(\pi\times T)(F)\|: (\pi,T)\ \textup{is a r-cov.con.
 pair}\}.
\]
\end{definition}
The left semicrossed product is endowed with an operator space
structure (the maximal one, see \cite[1.2.22]{BleLeM04}) induced by
the matrix norms
\[
 \|[F_{i,j}]\|_l=\sup\{\|[(T\times\pi)(F_{i,j})]\|: (\pi,T)\
 \textup{l-con.cov. pair}\}.
\]
 We note that there is a bijective correspondence between the
l-cov.con. pairs $(\pi,T)$ and the unital completely contractive
representations of $\ell^1(\bbZ_+,\C,\ga)_l$. Thus, the left
semicrossed product has the following universal property (up to
completely isometric isomorphisms): for any unital operator algebra
$\B$ and for any unital completely contractive morphism
$\rho:\ell^1(\bbZ_+,\C,\ga)_l\rightarrow\B$, there exists a unique
unital completely contractive morphism
$\tilde{\rho}:\bbZ_+\times_\ga\C\rightarrow\B$ that extends
$\rho$.\\

In theorem \ref{onlyone}, we prove that the semicrossed product, as
an operator algebra, is independent of the way $\C$ is (faithfully)
represented. In order to do so, we use some dilations theorems of
\cite{MuhSol06} and \cite{Pet84} and arguments similar to the ones
in \cite[theorem 6.2]{Kats04}.

First of all, every l-cov.con. pair $(\pi,T)$ on a Hilbert space $H$
dilates to a l-cov.is. pair $(\eta,W)$ on a Hilbert space
$H_1\supseteq H$, such that $\eta(x)H\subseteq H$ and
$\eta(x)|_H=\pi(x)$, for every $x\in\C$, and $T^n=P_HW^n|_H$, for
every $n\in\bbZ_+$, where $W$ is an isometry (see \cite{MuhSol06}).
Hence, by \cite[II.5]{Pet84} we see that the norm $\|\cdot\|_l$ is
the supremum over all left covariant purely isometric pairs. By
\cite[proposition I.4]{Pet84}, for such a pair $(\eta,W)$ on a
Hilbert space $H_1$ there is a representation $(H_2,\pi')$ of $\C$
such that $W\times\eta$ is unitarily equivalent to $lt_{\pi'}$.
Thus, eventually we have that, for $F\in \ell^1(\bbZ_+,\C,\ga)_l$,
$\|F\|_l=\sup\{\|lt_\pi(F)\|: (H,\pi) \text{ a *-representation of }
\C\}$. Moreover, $\|[F_{i,j}]\|_l=\sup\{\|[lt_\pi(F_{i,j})]\|:
(H,\pi) \text{ a *-representation of } \C\}$.

\begin{proposition}
If $F_{i,j}\in\ell^1(\bbZ_+,\C,\ga)_l$, then
$\|[F_{i,j}]\|_l=\|[lt_{\pi_u}(F_{i,j})]\|$, where $(H_u,\pi_u)$ is
the universal representation of $\C$. Analogously, for every
$F_{i,j}\in \ell^1(\bbZ_+,\C,\ga)_r$,
$\|[F_{i,j}]\|_r=\|[rt_{\pi_u}(F_{i,j})]\|$.
\end{proposition}
\proof Let $(H,\pi)$ be a *-representation of $\C$. By definition of
the universal representation we have that $\pi_u|_H=\pi$ and
$\pi_u(x)H\subseteq H$. Let $H_0=H\otimes\ell^2(\bbZ_+)$. We denote
by $P_{H_0}$ the projection onto $H\otimes\ell^2(\bbZ_+)\subseteq
H_u\otimes\ell^2(\bbZ_+)$ and observe that $P_{H_0}
(\mathbf{1}_{H_u}\otimes v)^n|_{H_0}= (\mathbf{1}_{H_0}\otimes
v)^n$, for every $n\in\bbZ_+$. Thus, for every $\nu\in\bbZ_+$ and
for every $[F_{i,j}]\in \M_\nu\left(\ell^1(\bbZ_+,\C,\ga)\right)$,
we have that $[lt_\pi(F_{i,j})]=(P_{H_0}\otimes I_\nu)
[lt_{\pi_u}(F_{i,j})]|_{(H_0)^{(\nu)}}$, and so
$\|[lt_\pi(F_{i,j})]\|\leq\|[lt_{\pi_u}(F_{i,j})]\|$. $\Box$\\

If $(H,\pi)$ is a faithful *-representation of $\C$, we denote by
$C^*(\pi,V)$ the C*-algebra generated by the representation $lt_\pi$
in $\B(H\otimes\ell^2(\bbZ_+))$. The covariance relation shows that
$C^*(\pi,V)$ is the norm-closed linear span of the monomials
$V^m\widetilde{\pi}(x)(V^*)^\gl$, $m,\gl\in\bbZ_+$. Since,
$C^*(\pi,V)$ is a direct summand of $C^*(\pi_u,V_u)$, the
compression $\Phi:\B(H_u\otimes\ell^2(\bbZ_+))\rightarrow
\B(H\otimes\ell^2(\bbZ_+))$ is a *-epimorphism when restricted on
$C^*(\pi_u,V_u)$. We will prove that it is also faithful, hence
completely isometric.

To this end, for every $s\in[0,2\pi]$, we define
$u_s:\ell^2(\bbZ_+)\rightarrow\ell^2(\bbZ_+)$ by $u_s(e_m)=e^{2\pi
is}e_m$. Let $\widetilde{U_s}=\mathbf{1}_{H_u}\otimes u_s$ and
$U_s=\mathbf{1}_{H}\otimes u_s$. The map
$\widetilde{\gamma_s}=ad_{\widetilde{U_s}}$ is a *-automorphism of
$C^*(\pi_u,V_u)$, since
$\widetilde{\gamma_s}(\widetilde{\pi_u}(x))=\widetilde{\pi_u}(x)$
and $\widetilde{\gamma_s}(\widetilde{V_u}^n)=e^{2\pi
ins}\widetilde{V_u}^n$. Similarly, $\gamma_s=ad_{U_s}$ is a *-automorphism
of $C^*(\pi,V)$. It is clear that $\Phi\circ\widetilde{\gamma_s}=
\gamma_s\circ\Phi$, because $\Phi(\widetilde{U_s})=U_s$. We denote
by $C^*(\pi_u,V_u)^{\widetilde{\gamma}}$ the fixed point algebra of
$\widetilde{\gamma}$ and define the contractive, faithful projection
$\widetilde{E}:C^*(\pi_u,V_u)\rightarrow
C^*(\pi_u,V_u)^{\widetilde{\gamma}}$ by
  \[\widetilde{E}(X):=\int_0^{2\pi}
  \widetilde{\gamma_s}(X)\frac{ds}{2\pi},\]
(as a Riemann integral of a norm-continuous function). Let
$\B_k:=\{\sum_{n=0}^k
V_u^n\widetilde{\pi_u}(x_n)(V^*_u)^n:x_n\in\C\}$; then we can check
that $C^*(\pi_u,V_u)^{\widetilde{\gamma}}$ is the norm-closure of
$\cup_{k\in\bbZ_+}\B_k$. Let $X_k$ be an element of $\B_k$. Since,
$V_u^n\widetilde{\pi_u}(x)(V^*_u)^n= diag
\{\underbrace{0,\dots,0}_{\text{n-times}},\pi_u(x),
\pi_u(\ga(x)),\dots\}$, we see that $X_k$ is a diagonal matrix whose
$(m,m)$-entry is the element $({X_k})_{m,m}= \pi_u \left(
\sum_{j=0}^{\min\{m,k\}} \ga^{m-j}(x_{m-j}) \right)$. So, if
$(H,\pi)$ is a faithful *-representation of $\C$,
\begin{align*}
\|(X_k)_{m,m}\| & =  \|\pi_u \left( \sum_{j=0}^{\min\{m,k\}}
\ga^{m-j}(x_{m-j}) \right)\|\\ &= \|\sum_{j=0}^{\min\{m,k\}}
\ga^{m-j}(x_{m-j})\|_\C= \|\pi \left( \sum_{j=0}^{\min\{m,k\}}
\ga^{m-j}(x_{m-j}) \right)\| \\ &= \|(\Phi(X_k))_{m,m}\|.
\end{align*}
So $
\|X_k\|=\sup_m\{\|(X_k)_{m,m}\|\}=\sup_m\{\|(\Phi(X_k))_{m,m}\|\}=
\|\Phi(X_k)\|$; hence $\Phi:C^*(\pi_u,V_u)\rightarrow C^*(\pi,V)$ is
isometric on each $\B_k$. Thus, \emph{$\Phi$ is injective when
restricted to the fixed point algebra
$C^*(\pi_u,V_u)^{\widetilde{\gamma}}$}.

\begin{theorem}
\label{onlyone} The left semicrossed product $\bbZ_+\times_\ga\C$ is
completely isometrically isomorphic to the norm-closed linear span
of $\sum_{n=0}^k V^n\widetilde{\pi}(x_n)$, $x_n\in \C$, where
$(H,\pi)$ is any faithful *-representation of $\C$. Respectively,
the right semicrossed product $\C\times_\ga\bbZ_+$ is completely
isometrically isomorphic to the norm-closed linear span of
$\sum_{n=0}^k \widetilde{\pi}(x_n)(V^*)^n$, $x_n\in \C$, where
$(H,\pi)$ is any faithful *-representation of $\C$.
\end{theorem}
\proof It suffices to prove that the natural *-epimorphism $\Phi$ is
faithful, hence a (completely) *-isometric isomorphism. Let
$X\in\ker\Phi$, then $X^*X\in\ker\Phi$. Hence,
\begin{align*}
\Phi(\widetilde{E}(X^*X))&=\Phi\left(\int_0^{2\pi}
  \widetilde{\gamma_s}(X^*X)\frac{ds}{2\pi}\right)\\
  &=\int_0^{2\pi}
  \Phi(\widetilde{\gamma_s}(X^*X))\frac{ds}{2\pi}
  = \int_0^{2\pi}\gamma_s(\Phi(X^*X))\frac{ds}{2\pi}=0.
\end{align*}
Now $\widetilde{E}(X^*X)$ is in
$C^*(\pi_u,V_u)^{\widetilde{\gamma}}$ and $\Phi$ is faithful there;
hence $\widetilde{E}(X^*X)=0$ and so
$X^*X=0$.\\
For the right semicrossed product, note that
$C^*(\pi,V^*)=C^*(\pi,V)$. $\Box$\\

If, in particular, $\ga$ is a *-isomorphism, then there is a natural
way to identify the left semicrossed product as a closed subalgebra
of the (reduced) crossed product, i.e. $C^*_l(\C)$. In this case, we
refer to this closed subalgebra as the \emph{left reduced
semicrossed product}. In a dual way, we can define the \emph{right
reduced semicrossed product}. The following is proved in
\cite{Petarx}, when $\C$ is abelian.

\begin{theorem}
If $\ga$ is a *-isomorphism, then the C*-envelope of the semicrossed
product is the (reduced) crossed product.
\end{theorem}
\proof Since $\ga$ is a *-isomorphism, we can view
$\ell^1(\bbZ_+,\C,\ga)_l$ as a $|\cdot|_1$-closed subalgebra of
$\ell^1(\bbZ,\C,\ga)_l$. First we prove that the inclusion map
$\ell^1(\bbZ_+,\C,\ga)\hookrightarrow\ell^1(\bbZ,\C,\ga)$ is
completely isometric. The key is to prove that
\[
 \|F\|_l=\sup\{\|(U\times\pi)(F)\|: (\pi,U) \text{ l-cov.un. pair of }
  \ell^1(\bbZ,\C,\ga)_l\},
\]
for every $F\in\ell^1(\bbZ_+,\C,\ga)_l$, since the right hand side
is exactly the norm of the (left) crossed product. For simplicity,
we denote this norm by $\|\cdot\|$. It is obvious that $\|F\|\leq
\|F\|_l$, since every l-cov.un pair of $\ell^1(\bbZ,\C,\ga)_l$
restricts to a l-cov.un. pair of the subalgebra
$\ell^1(\bbZ_+,\C,\ga)_l$. Also, if $(H_0,\pi)$ is a faithful
*-representation of $\C$, then $lt_\pi$ is the compression of the
left regular representation of $\ell^1(\bbZ,\C,\ga)_l$ induced by
$\pi$, denoted simply by $lt$. So, $\|lt_\pi(F)\|\leq\|lt(F)\|$,
thus $\|F\|_l\leq\|F\|$ by theorem \ref{onlyone}. Arguing in the
same way, we get that $\|[F_{i,j}]\|\leq\|[F_{i,j}]\|_l$ and
$\|[lt_\pi(F_{i,j})]\|$ $\leq \|[lt(F_{i,j})]\|$, for every
$[F_{i,j}]\in \M_\nu \left(\ell^1(\bbZ_+,\C,\ga)_l\right)$. But $lt$
is a *-morphism of the crossed product, hence completely
contractive. Thus, $\|[F_{i,j}]\|_l\leq\|[F_{i,j}]\|$ and equality
holds.\\
Hence, if $\widehat{\pi}(x)=diag\{\pi (a^m(x)),m\in\bbZ\}$ and
$U=1_{H_0}\otimes u$, where $u$ is the bilateral shift, then the map
$\gd_n\otimes x\mapsto U^n\widehat{\pi}(x)$ extends to a complete
isometry $\gi:\bbZ_+\times_\ga\C\rightarrow C^*_l(\C)$, whose image
generates $C^*_l(\C)$ as a C*-algebra. Let $\B$ be the C*-envelope
of $\bbZ_+\times_\ga\C$. Then, by the universal property of
C*-envelopes, there is a surjective C*-homomorphism $\Psi: C^*_l(\C)
\rightarrow \B$, which restricts to a completely isometry on
$\gi(\bbZ_+\times_\ga\C)$. Let $G\in \ker\Phi$ be of unit norm, and
choose $F=\sum_{n=-k}^k U^n\widehat{\pi}(x_n)$ with $\|G-F\|<
\frac{1}{2}$. Thus $U^kG\in\ker\Psi$,
$\gi^{-1}(U^kF)\in\bbZ_+\times_\ga\C$,
$\|\gi^{-1}(U^kF)\|=\|U^kF\|=\|F\|>\frac{1}{2}$ and
$\|U^kG-U^kF\|=\|G-F\|<\frac{1}{2}$. Then $
 \frac{1}{2}<\|i^{-1}(U^kF)\|=\|\Psi(U^kF)\|=\|\Psi(U^kF-U^kG)\|\leq
 \|U^kF-U^kG\|<\frac{1}{2},
$ which is a contradiction. $\Box$

\section{w*-Semicrossed products}

\noindent

Let $\A\subseteq B(H_o)$ be a unital subalgebra, closed in the
w*-operator topology, and $\beta:\A\rightarrow\A$, a contractive
w*-continuous endomorphism of $\A$. From now on we fix $H=H_o\otimes
\ell ^2(\bbZ_+)$ and $\pi:=\widetilde{id_\A}$, as in example 1.1.
Then $\pi $ is a faithful representation of $\A$ on $H$, and we can
write $\pi (b)=\sum_{n\geq0} \beta ^n (b)\otimes p_n$, where $p_n $
$\in \B(\ell ^2(\bbZ_+))$ is the projection onto $[e _n]$. Note that
the sum converges in the w*-topology as well. Hence, $\pi (b)$
belongs to the \emph{w*-tensor product algebra} $\A\overline{\otimes
}\B(\ell ^2(\bbZ_+))$. This is, by definition, the w*-closed linear
span in $\B(H)$ of the operators $b\otimes a$, with $b\in \A$ and
$a\in \B(\ell ^2(\bbZ_+))$. We also represent $\bbZ_+$ on $H$ by the
isometries $V^n= \mathbf{1}_{H_o}\otimes v^n$, where $v$ is the
unilateral shift on $\ell ^2(\bbZ_+)$. Thus, $V^n\in
\A\overline{\otimes }\B(\ell ^2(\bbZ_+))$.

\begin{definition}
The w*-semicrossed product $\bbZ_+\overline{\times}_\beta \A $ is
the w*-closure of the linear space of the `analytic polynomials' $
\sum_{n=0}^k V^n\pi(b_n)$, $b_n\in \A$, $k\geq 0$.
\end{definition}

It is easy to check that the left covariance relation $
 \pi(b)V=V\pi(\beta (b))$ holds. Hence, $(\pi ,V)$ is a left covariant
isometric pair. Thus, the w*-semicrossed product is a unital
(non-selfadjoint) subalgebra of $\B(H)$ and by definition, $\bbZ_+
\overline{\times}_\beta \A \subseteq \A\overline{\otimes }\B(\ell
^2(\bbZ_+))$.

\begin{proposition}
\label{char}An operator $T\in \B(H)$ is in the w*-semicrossed
product if and only if $T_{\kappa,\gl}\in\A$ and
$G_m(T)=V^m\pi(T_{m,0})$ when $m\in \bbZ_+$ while $G_m(T)=0$ for
$m<0$. Equivalently, when $T_{\kappa,\gl}\in\A$ and
$\beta(T_{m+\gl,\gl})=T_{m+\gl+1,\gl+1}$ for every $m,\gl\in\bbZ_+$,
while $T_{\kappa,\gl}=0$ when $\kappa<\gl$.
\end{proposition}
\proof If $T=\sum_{\kappa=0}^nV^\kappa\pi(b_\kappa)$ with
$b_\kappa\in \A$, then $G_m(T)=V^m\pi(b_m)$ when
$m\in\{0,1,\dots,n\}$ and $G_m(T)=0$ otherwise. Let
$T\in\bbZ_+\overline{\times}_\beta \A$ and a net
$A_i=\sum_{\kappa=0}^{n_i}V^\kappa\pi(b_{i,\kappa})$ of analytic
polynomials converging to $T$ in the w*-topology. Since $G_m$ is
w*-continuous, we have that $G_m(T)=\text{w*-}\lim_i G_m(A_i)$ for
every $m\in\bbZ$. Thus $G_m(T)=0$ when $m<0$. If $m\geq0$, then
$T_{(m)}=(V^*)^mG_m(T)=\text{w*-}\lim_i(V^*)^mG_m(A_i)=\text{w*-}
\lim_i\pi(b_{i,m})$. Let $\phi\in\B(H_0)_*$ and $k\in\bbZ_+$, then
$\phi\overline{\otimes}\om_{e_\kappa,e_\kappa}\in\B(H)_*$; hence we
get $\phi(T_{m+\kappa,\kappa})=
(\phi\overline{\otimes}\om_{e_\kappa,e_\kappa}) \left(T_{(m)}
\right)=\lim_i
(\phi\overline{\otimes}\om_{e_\kappa,e_\kappa})(\pi(b_{i,n}))=
\lim_i\phi(\beta^k(b_{i,n})).$ Thus
$T_{m+\kappa,\kappa}=\text{w*-}\lim_i \beta^\kappa(b_{i,m})$, for
every $\kappa\in\bbZ_+$, so $T_{m+\kappa,\kappa}\in \A$. Also, since
$\beta$ is w*-continuous, we get that
$\beta^\kappa(T_{m,o})=\text{w*-}\lim_i \beta^\kappa(b_{i,m})=
T_{m+\kappa,\kappa}$, for every $\kappa\in\bbZ_+$. Hence, we get
that $G_m(T)=V^m
\pi(T_{m,o})$, for every $m\geq0$.\\
For the opposite direction, if $T\in\B(H)$ satisfies the conditions,
we can see that $G_m(T)\in\bbZ_+\overline{\times}_\beta \A$, and so
by the F\'{e}jer Lemma, $T\in\bbZ_+\overline{\times}_\beta \A$ as
well. The last equivalence is trivial. $\Box$

\begin{remark}
\label{rem}\textup{Note that each $ad_{U_s}$ leaves
$\bbZ_+\overline{\times}_\beta \A$ invariant, and hence, being
unitarily implemented, also leaves its reflexive cover invariant.
Thus, so does $G_m(\cdot)$.}
\end{remark}

Suppose now that the endomorphism $\beta$ is \emph{implemented }by a
unitary $w$ acting on $H_0$, so that $ \beta (b)=w bw^*$, for all
$b\in \A$. Let $\rho (b)=b\otimes \mathbf{1}_{\ell ^2(\bbZ_+ )}$,
for $b\in\A$ and $W =w^* \otimes v$. Then $(\rho ,W)$ is a left
covariant isometric pair and we denote by $\bbZ_+\times_w \A$
\emph{the w*-closure of the linear space of the `analytic
polynomials'} $\sum_{n=0}^k W^n\rho(b_n)$, $b_n\in\A$, $k\geq 0$.

It is easy to check that $\bbZ_+\times_w \A$ is unitarily equivalent
to $\bbZ_+\overline{\times}_\beta \A $, via $Q=\sum_{n\geq0}w^{-n
}\otimes p_n,$. Thus we refer to $\bbZ_+\overline{\times}_w \A$ as
the w*-semicrossed product, as well. Using the unitary operator $Q$
and proposition \ref{char} we get the following characterization.

\begin{proposition}
\label{char2}An operator $T\in \B(H)$ is in
$\bbZ_+\overline{\times}_w \A$ if and only if $G_m(T)=W^m\rho(b_m)$,
for some $b_m\in\A$, when $m\in \bbZ_+$ and $G_m(T)=0$ for $m<0$.
Equivalently, when $T_{m+\gl,\gl}=(w^*)^mb_m$, for every
$m,\gl\in\bbZ_+$ and $T_{\kappa,\gl}=0$, when $\kappa<\gl$. $\Box$
\end{proposition}
The relation between the w*-tensor product $\A\overline{\otimes} \T$
and $\bbZ_+\overline{\times}_w \A$ depends on some properties of
$w$. Specifically,
\begin{itemize}
\item $\A\overline{\otimes}\T=\bbZ_+\overline{\times}_w \A$ if and only if
 $w, w^*\in\A$.
\item $\A\overline{\otimes}\T\varsubsetneq \bbZ_+\overline{\times}_w \A$
 if and only if $w^*\notin\A$, $w\in\A$.
\item $\bbZ_+\overline{\times}_w \A\varsubsetneq \A\overline{\otimes}\T$
 if and only if $w\notin\A$, $w^*\in\A$.
\item $(\bbZ_+\overline{\times}_w \A)\cap(\A\overline{\otimes}\T)=\rho(\A)$,
 if and only if $(w^n\A)\cap\A=\{0\}$, $\forall \ n\in\bbZ_+$.
\end{itemize}
It is easy to verify that, when $(w^n\A)\cap\A=\{0\}$ for every
$n\in\bbZ_+$, then $w,w^*\notin\A$, but the converse is not always
true.

\begin{example}
\textup{Take $\A=L^\infty(\bbT)$ acting on $L^2(\bbT)$ and
$\beta(f)(z)=f(\gl z)$, where $\gl$ is a q-th root of unity. Then
$\beta$ is unitarily implemented by $w\in\B(L^2(\bbT))$, with
$(w(g))(z)=g(\gl z)$. Then $w^{mq}=I_{H_0}$, for every $m\in\bbZ_+$,
hence $w^{mq}\A\cap\A=\A$. In this case, $(\bbZ_+\overline{\times}_w
\A)\cap(\A\overline{\otimes}\T)$ contains the w*-closed algebra
generated by $\sum_{n=0}^k W^{nq}\rho(b_n)$, $ b_n\in\A$, which
properly contains $\rho(\A)$. }
\end{example}
The following lemma will be superseded below (theorem \ref{ref}).

\begin{lemma}
\label{reflemma} The w*-semicrossed product
$\bbZ_+\overline{\times}_w\B(H_0)$ is reflexive, for every unitary
$w\in\B(H_0)$.
\end{lemma}
\proof Let $T\in \Ref(\bbZ_+\overline{\times} _w \B(H_0))$. By
remark \ref{rem} each $G_m(T)$ belongs to the reflexive cover of the
w*-semicrossed product. For $\kappa<\gl$ and $\xi,\eta\in H_0$ there
is a sequence $A_n\in\bbZ_+\overline{\times}_w\B(H_0)$ such that
$\sca{T(\xi\otimes e_\gl),\eta\otimes
e_{\kappa}}=\lim_n\sca{A_n(\xi\otimes e_\gl),\eta\otimes
e_{\kappa}}$. Hence, $ \sca{T_{\kappa,\gl}\xi,\eta}
=\lim_n\sca{(A_n)_{\kappa,\gl}\xi,\eta}=0, $ since each
$(A_n)_{\kappa,\gl}=0$, for $\kappa<\gl$. Hence $G_m(T)=0$ for every
$m<0$. Now, fix $m\in\bbZ_+$ and consider $\xi\in H_0$,
$g_r=\sum_nr^ne_n$, $0\leq r<1$. We can check that the subspace
$\F=\overline{[(b\xi)\otimes g_r: b\in\B(H_0)]}$ is
$(\bbZ_+\overline{\times}_w\B(H_0))^*$-invariant, hence
$G_m(T)^*$-invariant. Since $\xi\otimes  g_r \in \F$, there is a
sequence $(b_n)$ in $\B(H_0)$ such that $G_m(T)^*(\xi\otimes
g_r)=\lim_n (b_n\xi)\otimes g_r$. Thus, $\sum_\kappa
r^{m+\kappa}T_{m+\kappa,\kappa}^*\xi\otimes e_\kappa= \lim_n
(b_n\xi)\otimes g_r$. Taking scalar product with $\eta\otimes
e_{\kappa}$, where $\eta\in H_0$ and $\kappa\geq 0$, we have that
$r^{m+\kappa}\sca{T_{m+\kappa,\kappa}^*\xi,\eta}=\lim_n
r^\kappa\sca{b_n\xi,\eta}$. Hence, $
r^{m}\sca{T_{m+\kappa,\kappa}^*\xi,\eta}=\lim_n
\sca{b_n\xi,\eta}=r^m\sca{T_{m,0}^*\xi, \eta}$, for every $\eta$.
Thus, $T_{m+\kappa,\kappa}^*\xi=T_{m,0}^*\xi$, for arbitrary $\xi\in
H_0$, so $T_{m+\kappa,\kappa}=T_{m,0}$ for every $\kappa\in\bbZ_+$.
Hence, $G_m(T)\in\B(H_0)\overline{\otimes} \T$, which coincides with
$\bbZ_+\overline{\times}_w\B(H_0)$
since $w\in\B(H_0)$. $\Box$\\

Let $\S$ be a w*-closed subspace of $\B(H)$. We say that $\S$ is
\emph{$G$-invariant} if $G_m(\S)\subseteq \S$ for every $m\in\bbZ$.
If, in particular, $\S$ is a w*-closed subspace of
$\bbZ_+\overline{\times}_w\B(H_0)$, then $G_m(\S)=0$, for every
$m<0$. In the next proposition we prove that we can associate a
sequence $(\S_m)_{m\geq0}$ of w*-closed subspaces of $\B(H_0)$ to
such an $\S$, and vice versa.

\begin{proposition}
A w*-closed subspace $\S$ of $\B(H)$ is a $G$-invariant subspace of
$\bbZ_+\overline{\times}_w\B(H_0)$ if and only if it is the
w*-closure of the linear space of the analytic polynomials
$\sum_{n=0}^kW^n\rho(x_n)$, $x_n\in \S_n, k\in\bbZ_+$, where $\S_n$
are w*-closed subspaces of $\B(H_0)$.
\end{proposition}
\proof Let $\S$ be a $G$-invariant w*-closed subspace of
$\bbZ_+\overline{\times}_w\B(H_0)$ and let $\S_m=\{w^mT_{m,0}:T\in
\S \}$, for every $m\geq 0$. Then $\S_m$ is a w*-closed subspace of
$\B(H_0)$. Indeed, let $x=\textup{w*-}\lim_i w^m(T_i)_{m,0}$, for
$T_i\in S$. Then $\rho((w^*)^mx)= \textup{w*-}\lim_i
\rho((T_i)_{m,0})$, so $W^m\rho(x)= \textup{w*-}\lim_i
V^m\rho((T_i)_{m,0})=\textup{w*-}\lim_i G_m(T_i)$, since
$T_i\in\bbZ_+\overline{\times}_w\B(H_0)$. But $\S$ is $G$-invariant,
hence $W^m\rho(x)\in \S$, thus $x=w^m(V^m\rho(x))_{m,0}\in \S_m$. A
use of the F\'{e}jer Lemma and proposition \ref{char2}, completes
the forward
implication.\\
For the converse, let $\S$ be a w*-closed subspace as in the
statement and $A\in \S$; so $A=\textup{w*-}\lim_i A_i$, where
$A_i=\sum_{\kappa=0}^{n_i}W^\kappa\rho(x_{i,\kappa})$, with
$x_{i,\kappa}\in \S_\kappa$. Then
$A\in\bbZ_+\overline{\times}_w\B(H_0)$ and
$G_m(A)=\textup{w*-}\lim_i G_m(A_i)=\textup{w*-}\lim_i
W^m\rho(x_{i,\kappa})$. Thus, $w^mA_{m,0}=\textup{w*-}\lim_i
x_{i,m}\in \S_m$. Hence, we have that $G_m(A)=W^m\rho(w^mA_{m,0})\in
\S$. $\Box$

\begin{theorem}
Let $(\S_m)_{m\geq 0}$ be the sequence associated to a $G$-invariant
w*-closed subspace $\S$ of $\bbZ_+\overline{\times}_w\B(H_0)$. If
every $\S_m$ is reflexive then $\S$ is reflexive.
\end{theorem}
\proof  By lemma \ref{reflemma}, $\Ref(\S)\subseteq
\Ref(\bbZ_+\overline{\times}_w\B(H_0))=\bbZ_+\overline{\times}_w\B(H_0)$.
So, for every $T$ in the reflexive cover of $\S$ and every
$m,\gl\in\bbZ_+$, we have that $T_{m+\gl,\gl}=(w^*)^mb_m$, where
$b_m\in\B(H_0)$. Thus, it suffices to prove that $b_m\in \S_m$.
Since $T\in \Ref(\S)$, for every $\xi,\eta\in H_0$, there is a
sequence $(A_n)$ in $\S$ such that $\sca{T(\xi\otimes
e_\gl),(w^*)^{m}\eta\otimes e_{m+\gl}}=\lim_n \sca{A_n(\xi\otimes
e_\gl),(w^*)^{m}\eta\otimes e_{m+\gl}}$. So,
$\sca{b_m\xi,\eta}=\sca{T_{m+\gl,\gl}\xi,(w^*)^{m}\eta}=\lim_n
\sca{(A_n)_{m+\gl,\gl}\xi,(w^*)^{m}\eta}$. Since each $A_n\in \S$,
then $(A_n)_{m+\gl,\gl}=(w^*)^{m}b_{n,m}$ for some $b_{n,m}\in
\S_m$. Thus $\sca{b_m\xi,\eta}=\lim_n\sca{b_{n,m}\xi,\eta}$, which
means that $b_m\in \Ref(\S_m)=\S_m$. $\Box$

\begin{theorem}
\label{ref}If $\A$ is a reflexive algebra, then
$\bbZ_+\overline{\times}_w \A$ is reflexive. In addition, if $\A$ is
hereditarily reflexive, then every $G$-invariant w*-closed subspace
of $\bbZ_+\overline{\times}_w \A$ is reflexive.
\end{theorem}
\proof The algebra $\bbZ_+\overline{\times}_w \A$ is associated to
the sequence $(\A)_{m\geq 0}$; hence it is reflexive by the previous
theorem. $\Box$

\begin{applications}
\textbf{A.} \textup{ (\textit{Sarason's result}, \cite[theorem
3]{Sar66}) Consider the case of a reflexive subalgebra $\A$ of
$M_n(\bbC)$ and a unitary $w\in M_n(\bbC)$ such that $w\A
w^*\subseteq \A$. Then $\bbZ_+\overline{\times}_w \A$ is reflexive.
Note that $\bbZ_+\overline{\times}_w \A=\T$ when $n=1$ and
$w=I_{H_0}$.}
\medskip

\textbf{B.} \textup{ (\textit{Ptak's result}, \cite[theorem
2]{Pta86}) More generally, $\A\overline{\otimes}\T$ coincides with
$\bbZ_+\overline{\times}_{I_{H_0}}\A$. So $\A\overline{\otimes}\T$
is reflexive, when $\A$ is reflexive.}
\medskip

\textbf{C.} \textup{ If $\M$ is a maximal abelian selfadjoint
algebra and let $\beta$ be a *-automorphism, then
$\bbZ_+\overline{\times}_\beta \M$ is reflexive, since every
*-automorphism of a m.a.s.a. is unitarily implemented. For example
let $\M=L^\infty(\bbT)$ acting on $L^2(\bbT)$ and $\beta$ the
rotation by $\theta\in\bbR$. Also $\bbZ_+\overline{\times}_\beta \A$
is reflexive whenever $\A$ is a $\beta$-invariant w*-closed
subalgebra of $\M$, since $\M$ is hereditarily reflexive (see
\cite{LogSu75}).}
\medskip

\textbf{D.} \textup{ Consider $\T$ acting on $H^2(\bbT)$ and $\beta$
as in the previous example. Then, $\T$ is reflexive and so
$\bbZ_+\overline{\times}_\beta \T$ is a reflexive subalgebra of
$\B(H^2(\bbT))\overline{\otimes}\B(\ell^2(\bbZ_+))$.}
\medskip

\textbf{E.} \textup{ If $\A$ is a nest algebra and $\beta$ is an
isometric automorphism, then it is unitarily implemented (see
\cite{Dav88}). Thus, $\bbZ_+\overline{\times}_\beta \A$ is
reflexive.}
\medskip

\textbf{F.} \textup{ Consider a C*-algebra $\C$ and a *-morphism
$\ga:\C\rightarrow\C$. Let $(H_0,\sigma)$ be a faithful
*-representation of $\C$ such that the induced *-morphism
\[
\beta:\sigma(\C)\rightarrow \sigma(\C):\sigma(x)\mapsto
\beta(\sigma(x))=\sigma(\ga(x))
\]
is implemented by a unitary $w\in\B(H_0)$. Then the induced
representation $lt_\sigma$ is faithful on $\bbZ_+\times_\ga\C$.
Thus, $\overline{lt_\sigma(\bbZ_+\times_\ga\C)}^{\textup{w*}}$ is
the w*-closed linear span of the analytic polynomials $\sum_{n=0}^k
V^n\pi(\sigma(x))$, and it is unitarily equivalent to the algebra
$\fC:= \overline{\text{span} \{\rho(\sigma(x)),W^n:x\in\C,
n\in\bbZ_+\}}^{\textup{w*}}$, via $Q=\sum_{n\geq0}w^{-n}\otimes
p_n$. But $\fC$ is exactly the w*-semicrossed product
$\bbZ_+\overline{\times}_w \overline{\sigma(\C)}^{\textup{w*}}$.
Thus, $\overline{lt_\sigma(\bbZ_+\times_\ga\C)}^{\textup{w*}}$ is
reflexive. \\
In particular, let $K$ be a compact, Hausdorff space,
$\mu$ a positive, regular Borel measure on $K$ and
$\sigma:C(K)\rightarrow \
B(L^2(K,\mu)):f\mapsto M_f$. Consider a
homeomorphism $\phi$ of $K$, such that $\phi$ and $\phi^{-1}$
preserve the $\mu$-null sets and let $\ga(f)=f\circ\phi$. Then the map
$M_f\rightarrow M_{f\circ\phi}$ extends to a *-automorphism of
$L^\infty(K,\mu)$, hence it is unitarily implemented.
Thus, $\overline{lt_\sigma(\bbZ_+\times_\ga C(K))}^{\textup{w*}}$ is
reflexive.}
\medskip

\textbf{H.} \textup{ Let $(\M,\tau)$ be a von Neumann algebra with a
faithful, normal, tracial state $\tau$ and let $L^2(\M,\tau)$ be the
Hilbert space associated to $(\M,\tau)$. Let $\beta:\M\rightarrow\M$
be a trace-preserving *-automorphism and consider $\M$ acting on
$L^2(\M,\tau)$ by left multiplication. Then $\beta$ is unitarily
implemented and it can be verified that the w*-semicrossed product
$\bbZ_+\overline{\times}_w \M$ coincides with the adjoint of the
analytic semicrossed product defined in \cite{McAsMuhSai79} and
\cite{Pel80}. Hence, we obtain \cite[proposition 4.5]{Pel80} for
$p=2$.}
\end{applications}

We conclude the analysis of the w*-semicrossed product
$\bbZ_+\overline{\times}_w \A$ by finding its commutant. We know
that $U_s(\bbZ_+\overline{\times}_w
\A)U_s^*=\bbZ_+\overline{\times}_w \A$, for all $s\in[0,2\pi]$,
hence, $U_s(\bbZ_+\overline{\times}_w
\A)'U_s^*=(\bbZ_+\overline{\times}_w \A)'$. Thus,
$T\in(\bbZ_+\overline{\times}_w \A)'$ if and only if
$G_m(T)\in(\bbZ_+\overline{\times}_w \A)'$, for every $m\in\bbZ$.
Now, recall that $w\A w^*\subseteq \A$, hence $w^*\A'w\subseteq\A'$.
So, we can define the w*-semicrossed product
$\bbZ_+\overline{\times}_\gamma \A'$, where $\gamma\equiv
ad_{w^*}:\A'\rightarrow\A'$.

\begin{theorem}
If $\gamma\equiv ad_{w^*}$, then $(\bbZ_+\overline{\times}_w
\A)'=\bbZ_+\overline{\times}_\gamma \A'$. $\Box$
\end{theorem}
\proof Obviously $T\in(\bbZ_+\overline{\times}_w \A)'$ if and only
if $T\in\{b\otimes\mathbf{1}, w^*\otimes v: b\in\A\}'$; note also
that $V\in(\bbZ_+\overline{\times}_w \A)'$. Let $T\in
(\bbZ_+\overline{\times}_w \A)'$, then for $m\geq 0$ and $b\in\A$,
$n\geq 0$,
\begin{align*}
G_m(T)(b\otimes\mathbf{1})=(b\otimes\mathbf{1})G_m(T)
&\text{ and } G_m(T)(w^*\otimes v)=(w^*\otimes v)G_m(T),\\
\text{hence, } T_{m+n,n}b=bT_{m+n,n} &\text{ and }
T_{m+n+1,n+1}(w^*)^n=(w^*)^nT_{m+n,n},\\
\text{so, } T_{m+n,n}\in\A' &\text{ and }
T_{m+n,n}=\gamma^n(T_{m,0}).
\end{align*}
Thus, if we set $\pi'(T_{m,0})=\sum_{n\geq0}\gamma^n(T_{m,0})\otimes
p_n$, we get that $G_m(T)=V^m\pi'(T_{m,0})$, for $m\geq 0$.\\
Now, let $m<0$, hence $G_m(T)= T_{(m)}(V^*)^{-m}$. Since,
$V\in(\bbZ_+\overline{\times}_w \A)'$, we have that $T_{(m)}
=G_m(T)V^{-m}\in(\bbZ_+\overline{\times}_w \A)'$. Thus,
$G_0(T_{(m)}) = T_{(m)} \in(\bbZ_+\overline{\times}_w \A)'$ and so,
by what we have proved, $T_{n,-m+n}=\gamma^n(T_{0,-m})$. Since
$G_m(T)(\xi\otimes e_0)=0$, then
$$G_m(T)(w^*\otimes
v)^{-m}(\xi\otimes e_0)= (w^*\otimes v)^{-m}G_m(T)(\xi\otimes
e_0)=0,$$
so $(T_{0,-m}(w^*)^{-m}\xi)\otimes e_0=0$. Thus,
$T_{0,-m}=0$ and therefore $T_{n,-m+n}=\gamma^n(T_{0,-m})=0$, for
every $n\geq 0$; hence $G_m(T)=0$, for every $m<0$. Hence, by
proposition \ref{char}, we get
that $T\in\bbZ_+\overline{\times}_\gamma \A'$.\\
For the converse, let $T\in\bbZ_+\overline{\times}_\gamma \A'$, then
$G_m(T)\in\bbZ_+\overline{\times}_\gamma \A'$ for every $m\in\bbZ$,
and we can see that $G_m(T)\in\{b\otimes\mathbf{1},w^*\otimes v:
b\in \A\}'$. Hence, $G_m(T)\in(\bbZ_+\overline{\times}_w \A)'$ for
every $m\in\bbZ$, so $T\in(\bbZ_+\overline{\times}_w \A)'$. $\Box$

\begin{theorem}
The double commutant of $\bbZ_+\overline{\times}_w \A$ is
$\bbZ_+\overline{\times}_w \A''$. Thus, the w*-semicrossed product
is its own bicommutant if and only if $\A=\A''$.
\end{theorem}
\proof We recall that $Q(\bbZ_+\overline{\times}_\beta \A)Q^*=
\bbZ_+\overline{\times}_w \A$, where $Q=\sum_{n}w^{-n}\otimes p_n$;
hence $Q^*(\bbZ_+\overline{\times}_\gamma \A')Q=
\bbZ_+\overline{\times}_{w^*} \A'$. Thus,
$(\bbZ_+\overline{\times}_w \A)''= (\bbZ_+\overline{\times}_\gamma
\A')'= (Q(\bbZ_+\overline{\times}_{w^*} \A')
Q^*)'=Q(\bbZ_+\overline{\times}_{w*} \A')'Q^*=
Q(\bbZ_+\overline{\times}_{\beta} \A'')Q^*=
\bbZ_+\overline{\times}_w \A''$. $\Box$\\

We end this section with a note on the \emph{reduced w*-semicrossed
products} (see the definition below). Let $\M$ be a von Neumann
algebra acting on a Hilbert space $H_0$, $\beta$ a *-automorphism of
$\M$ and consider $\bbZ\overline{\rtimes}_\beta \M$ to be the usual
w*-crossed product, a von Neumann subalgebra of
$\M\overline{\otimes}\B(\ell^2(\bbZ))$. This is by definition the
von Neuman algebra $\{\widehat{\pi}(b),U: b\in\M\}''$, where
$\widehat{\pi}(b)=\sum_{n\in\bbZ}\beta^n(b)\otimes p_n$ and
$U=1_{H_0}\otimes u$, the ampliation of the bilateral shift
$u\in\B(\ell^2(\bbZ))$.

\begin{definition}
The reduced w*-semicrossed product $\bbZ_+\overline{\rtimes} _\beta
\M$ is the w*-closure of the linear space of `analytic polynomials'
$\sum_{n=0}^k U^n\widehat{\pi}(b_n)$, $b_n\in \M$, $k\geq 0$.
\end{definition}

Since $(\widehat{\pi}, U)$ is a l-cov.un. pair, the reduced
w*-semicrossed product is a w*-closed subalgebra of the w*-crossed
product. In fact, note that $\bbZ_+\overline{\rtimes} _\beta \M$ is
the intersection of $\bbZ \overline{\rtimes} _\beta \M$ with the
`lower triangular' matrices. Hence, we have the following
proposition.

\begin{proposition}
The reduced w*-semicrossed product of a von Neumann algebra is
reflexive. $\Box$
\end{proposition}

Now, take $\A$ to be a w*-closed subalgebra of $\M$ which is
invariant under $\beta$. We define
$\bbZ_+\overline{\rtimes}_\beta\A$ to be \emph{the w*-closure of the
linear space of `analytic polynomials'} $\sum_{n=0}^k
U^n\widehat{\pi}(b_n)$, $b_n\in \A$, $k\geq 0$. Using the technique
of theorem \ref{ref} one can show the following.

\begin{corollary}
If $\A$ is reflexive subalgebra of $\M$ which is invariant under
$\beta$, then $\bbZ_+\overline{\rtimes}_\beta\A$ is reflexive.
$\Box$
\end{corollary}

\section{The commutative case}

\noindent

Now, we examine the case where $\C$ is a commutative, unital
C*-algebra, $\C=C(K)$, and the *-endomorphism $\ga$ is induced by a
continuous map $\phi:K\rightarrow K$. Let $ev_t$ be the evaluation
at $t\in K$, i.e. $ev_t(f)=f(t)$; then $(\ell^2(K),\oplus_tev_t)$ is
a faithful *-representation of $C(K)$. If some $t\in K$ has dense
orbit, we obtain a faithful representation of $C(K)$ on
$\ell^2(\bbZ_+)$. As observed in theorem \ref{onlyone}, such
representations play a fundamental role for the semicrossed product
$\bbZ_+\times_\ga C(K)$, since they are ``enough'' to obtain the
norm. Let $\pi_t:=\widetilde{ev_t}$, as in example 1.1. So,
$\pi_t:C(K)\rightarrow \B(\ell^2(\bbZ_+))$ is given by
$\pi_t(f):=\sum_{n\geq 0} f(\phi^n(t))p_n$, where $p_n$ is the
one-dimensional projection on $[e_n]$. Then $(\pi_t,v)$ is a left
covariant isometric pair. We define the \emph{one point
w*-semicrossed product} to be $\C_t=
\overline{lt_{\pi_t}(\bbZ_+\times_\ga C(K))}^{\textup{w*}}$, i.e.
the w*-closed linear span in $\B(\ell^2(\bbZ_+))$ of the `analytic
polynomials' $\sum_{n=0}^k
v^n\pi_t(f_n)$, $f_n\in C(K)$.\\

Let $t'=\phi^{n_0}(t)$ be the first periodic element of the orbit of
$t$ with period $p$, as in the following diagram
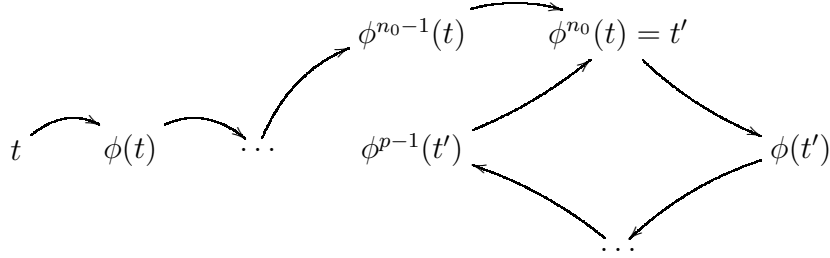
\[
\xymatrix{  & & & \phi^{n_0-1}(t)
\ar@/^1pc/[r] & \phi^{n_0}(t)=t' \ar@/_/[dr] & \\
t \ar@/^1pc/[r] & \phi(t) \ar@/^1pc/[r] &  \dots  \ar@/^1pc/[ur]
 & \phi^{p-1}(t') \ar@/_/[ur] & & \phi(t') \ar@/_/[dl] \\
& & & &  \dots \ar@/_/[ul] &}
\]
Then $orb(t)=\{t,\dots,\phi^{n_0-1}(t),t',\dots,\phi^{p-1}(t')\}$
induces a family of projections $\{P_{n_0},P_0,\dots,P_{p-1}\}$ such
that $I=P_{n_0}\oplus P_0 \oplus\cdots\oplus P_{p-1}$. Indeed, let
$P_{n_0}$ be the projection on $[e_0,\dots,e_{n_0-1}]$ and $P_i$ be
the projection on $[e_{n_0+i+pj}: j\in\bbZ_+]$ for $i=0,\dots,p-1$.
Note that if $f\in C(K)$, then $\pi_t(f)(e_{n_0+i+pj})=
f(\phi^{n_0+i+pj}(t))e_{n_0+i+pj} =f(\phi^{i}(t'))e_{n_0+i+pj}$, for
$j\in\bbZ_+$. Hence, $\pi_t(f)P_i= f(\phi^i(t'))P_i$, for every
$i=0,\dots p-1$.

\begin{proposition}
The algebra $\C_t$ is the linear sum $(\fT P_{n_0})\oplus (\T P_0)
\oplus \cdots \oplus (\T P_{p-1})$, where $\fT$ is the algebra of
lower triangular operators in $\B(\ell^2(\bbZ_+))$, $\T$ is the
algebra of analytic Toeplitz operators and $P_{n_0}, P_0,$ $\dots,
P_{p-1}$ are the projections induced by the orbit of $t$.
\end{proposition}
\proof For any $n\in\bbZ_+$ and $f\in C(K)$, we have
$$v^n\pi_t(f)= v^n\pi_t(f)P_{n_0}\oplus f(t')v^nP_0 \oplus
f(\phi^{p-1}(t'))v^nP_{p-1}.$$ Thus, $\C_t \subseteq (\fT
P_{n_0})\oplus (\T P_0) \oplus \cdots \oplus (\T P_{p-1})$.\\
For the converse, first let $TP_{n_0}\in  \fT P_{n_0}$ and note that
$(TP_{n_0})_{\kappa,\gl}=0$ when $\kappa<\gl$ or $n_0-1<\gl$. So,
$G_m(TP_{n_0})=0$, when $m<0$, and $G_m(TP_{n_0})=v^m
(\sum_{n=0}^{n_0-1} (TP_{n_0})_{m+n,n}p_n)$, when $m\geq 0$. Note
that $(TP_{n_0})_{\kappa,\gl}\in\bbC$ for every $\kappa,\gl\in
\bbZ_+$. Fix $m\geq 0$ and let $n\in\{0,\dots,n_0-1\}$. Then by
Urysohn's Lemma there is a sequence $(f_{n,j})_j$ of continuous
functions on $K$, such that $\lim_j
f_{n,j}(\phi^n(t))=(TP_{n_0})_{m+n,n}$ and $f_{n,j}(s)=0$ for $s\in
orb(t)\setminus\{\phi^n(t)\}$. Hence, $(TP_{n_0})_{m+n,n}p_n=
\text{w*-}\lim_j \pi_t(f_{n,j})\in \C_t$ and so
$v^m(TP_{n_0})_{m+n,n}p_n\in \C_t$. Thus $G_m(TP_{n_0})\in\C_t$,
and, by the F\'{e}jer Lemma, $TP_{n_0}\in \C_t$. So, $\fT
P_{n_0}\subseteq \C_t$.\\
Also, for fixed $i\in\{0,\dots,p-1\}$ and $m\in\bbZ_+$, consider
$v^mP_i\in \T P_i$. Again by Urysohn's Lemma, there is a sequence
$(f_{i,j})_j$ of continuous functions on $K$, such that $\lim_j
f_{i,j}(\phi^i(t'))=1$ and $f_{i,j}(s)=0$ for $s\in
orb(t)\setminus\{\phi^i(t')\}$. Then
$\text{w*-}\lim_j\pi_t(f_{i,j})=P_i$, so $v^mP_i\in\C_t$. Hence, $\T
P_i \subseteq \C_t$, for every $i\in\{0,\dots,p-1\}$. Thus, $(\fT
P_{n_0})\oplus (\T P_0) \oplus \cdots \oplus (\T
P_{p-1})\subseteq \C_t$. $\Box$\\
\\
Note that if $orb(t)$ has no periodic points, then $\C_t=\fT$, since
$P_{n_0}=\mathbf{1}_{\ell^2(\bbZ_+)}$. Also, if $orb(t)$ has exactly
one periodic point $t'$, then $\phi^n(t)=t'$ for every $n\geq n_0$
(i.e. $t'$ is a fixed point); thus $C_t=\fT P_{n_0}\oplus \T
P_{n_0}^\bot$. If $t$ is itself a fixed point, then $\C_t=\T$.

\begin{remark} \label{char3}
\textup{Let $\D$ be the algebra of diagonal operators
in $\B(\ell^2(\bbZ_+))$ and $\D_\phi= \{ T\in\D :
T_{\kappa,\kappa}=T_{n,n}\ \text{when} \ \phi^\kappa(t)=\phi^n(t)
\}$ which is a w*-closed subalgebra of $\D$. Hence, $T\in \D_\phi$
if and only if $T$ is of the form
$$T=diag\{y_0,\dots,y_{n_0-1},y_{n_0},
\dots,y_{p-1},y_{n_0},\dots,y_{p-1},\dots\}.$$ It is immediate from
the previous proposition that $\C_t$ is generated by the unilateral
shift in $\B(\ell^2(\bbZ_+))$ and the diagonal matrices id
$\D_\phi$. Thus, an operator $T\in\B(\ell^2(\bbZ_+))$ is in $\C_t$
if and only if for every $m<0$, $G_m(T)=0$, and for every $m\geq 0$,
$G_m(T)=v^m\sum_n T_{m+n,n}p_n$ where
$T_{m+\kappa,\kappa}=T_{m+n,n}$, whenever
$\phi^{\kappa}(t)=\phi^{n}(t)$.}
\end{remark}

\begin{theorem}
The algebra $\C_t$ is reflexive.
\end{theorem}
\proof If $T\in \Ref(\C_t)$, then $G_m(T)\in \Ref(\C_t)$; thus
$G_m(T)=0$, for $m<0$. Let $g_r=\sum_{n\geq0}r^ne_n$, with $0\leq
r<1$, and $\F=\overline{[\pi_t(f)g_r: f\in C(K)]}$. Then $\F$ is
$(\C_t)^*$-invariant; thus $G_m(T)^*$-invariant, for $m\in\bbZ_+$.
So, there is a sequence of $f_j\in C(K)$ such that $G_m(T)^*g_r=
\lim_j \pi_t(f_j)g_r$. Hence $r^{m}\overline{T}_{m+n,n}=\lim_j
(f_j(\phi^n(t)))$, for every $n\in\bbZ_+$. Thus,
$T_{m+n,n}=T_{m+\kappa,\kappa}$, if $\phi^\kappa(t)=\phi^n(t)$. So,
by remark \ref{char3}, $T\in\C_t$. $\Box$

\begin{remark} \label{unif}
\textup{In order to construct $\C_t$}, it is sufficient to take
coefficients
from any uniform algebra $\fA$ on $K$.\\
\textup{Indeed, let $\fA$ be a norm closed subalgebra of $C(K)$
containing the constant functions which separates the points of $K$
and form the polynomials $\sum_{n=0}^k v^n\pi_t(f_n)$, $f_n\in\fA$.
By remark \ref{char3}, it suffices to prove that $\pi_t(ball(\fA))$
is w*-dense in $ball(\D_\phi)$. Fix $z\in\bbT$ and $n_0\in\bbZ_+$,
and take $T\in \D_\phi$, such that $T_{n_0,n_0}=z$ and $T_{n,n}=1$,
if $\phi^n(t)\neq \phi^{n_0}(t)$. Using the argument of the claim of
\cite[theorem 2.9]{Kats04-2} we can find a sequence of $(f_j)_j$ in
$ball(\fA)$ such that $\textup{w*-}\lim_j\pi_t(f_j)=T$. To complete
the proof, observe that products of elements of this form
approximate the unitaries in $\D_\phi$ in the w*-topology and that
the strong closure of $\pi_t(ball(\fA))$ is closed under
multiplication.}
\end{remark}

\noindent \emph{\textbf{Acknowledgements.}} I wish to give my
sincere thanks to A. Katavolos for his kind help and advice during
the preparation of this paper. I also wish to thank E. Katsoulis for
bringing remark \ref{unif} to my attention. Finally, I wish to thank
I.Sis. and T.o.Ol. for the support and inspiration.

\bibliographystyle{plain}

\bibliography{biblio}

\end{document}